\documentclass{amsart} 
\usepackage{amssymb, amsmath}
\usepackage[dvips]{graphicx}
\usepackage{amsfonts}
\usepackage{latexsym}
\usepackage{color}
\newtheorem{theorem}{Theorem}	
\newtheorem{lemma}{Lemma}[section]		
		
\newtheorem{proposition}{Proposition}		
			
\newtheorem{definition}{Definition}	

\title{
Aperture of plane curves 
}
\author{
Daisuke Kagatsume}
\address{Machida Hall, Machida city, 
Tokyo 194-8520, Japan}
\email{kagatsume-daisuke-mt@ynu.jp}
\author{Takashi Nishimura
}
\address{Research Group of Mathematical Sciences,  
Research Institute of Environment and Information Sciences,  
Yokohama National University, 
Yokohama 240-8501, Japan}
\email{nishimura-takashi-yx@ynu.jp}
\begin{document}
\date{}
\begin{abstract}
For any given $C^\infty$ immersion ${\bf r}: S^1\to \mathbb{R}^2$ 
such that the set 
$\mathcal{NS}_{{\bf r}}=\mathbb{R}^2-\cup_{s\in S^1}\left({\bf r}(s)+d{\bf r}_{s}(T_s(S^1))\right)$ is not 
empty, a simple geometric model of crystal growth is constructed.    
It is shown that  
our geometric model of crystal growth never formulates a polygon while it is 
growing.   Moreover, it is shown also that our model always dissolves to a point.   
\end{abstract}
\maketitle  
\section{Introduction}
\label{introduction}
Let ${\bf r}: S^1\to \mathbb{R}^2$ be a $C^\infty$ immersion 
such that the set 
$$
\mathbb{R}^2-\bigcup_{s\in S^1}\left({\bf r}(s)+d{\bf r}_s(T_s(S^1))\right)
\leqno{(1.1)}
\label{(1.1)}
$$
is not the empty set, where  $T_{{\bf r}(s)}\mathbb{R}^2$ is identified with $\mathbb{R}^2$.                 
The perspective projection of the given plane curve ${\bf r}(S^1)$ from any point of (1.1) 
does not give the silhouette of  
${\bf r}(S^1)$ because it is non-singular.   
By this reason, the set {(1.1)} is called the {\it no-silhouette} of ${\bf r}$ 
and is denoted by $\mathcal{NS}_{\bf r}$ (see Figure 1).     
\begin{figure}[hbtp]
\begin{center}
\includegraphics[width=4cm]{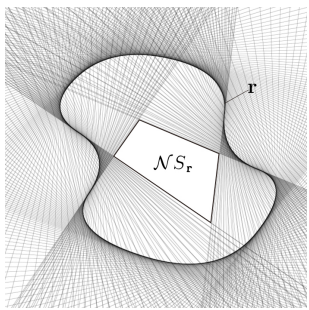}
\caption{The no-silhouette ${\mathcal NS}_{\bf r}$.}
\label{figure 1}
\end{center}
\end{figure}     
The notion of no-silhouette was first defined and studied from the viewpoint of 
perspective projection in \cite{nishimurasakemi1}.       
In \cite{nishimurasakemi2} 
it has been shown that the topological closure of no-silhouette is 
a Wulff shape,  which is the well-known geometric model of crystal at 
equilibrium introduced by G.~Wulff in \cite{wulff}.    
\par 
In this paper, we show that 
by rotating all tangent lines about their tangent points simultaneously with the same angle,   
we always obtain a geometric model of crystal growth (Proposition \ref{proposition 4}), 
our model never formulates a polygon while it is growing (Theorem \ref{theorem 1}), 
our model always dissolves to a point (Theorems \ref{theorem 2}), and our model is growing in a relatively simple way 
when the given ${\bf r}$ has no inflection points (Theorem 3).          
\par 
For any $C^\infty$ immersion ${\bf r}: S^1\to \mathbb{R}^2$ 
and any real number $\theta$, 
define the new set 
$$
\mathcal{NS}_{\theta, {\bf r}}=
\mathbb{R}^2-\bigcup_{s\in S^1}\left({\bf r}(s)+R_\theta \left(d{\bf r}_s(T_s(S^1))\right)\right),  
$$
where $R_\theta: \mathbb{R}^2\to \mathbb{R}^2$ is the rotation defined by 
$R_\theta(x,y)=(x\cos \theta -y\sin\theta, x\sin\theta +y\cos \theta)$ (see Figure 2).    
 \begin{figure}[htbp]
  \begin{center}
  \begin{minipage}{.45\linewidth} 
	\begin{center}
\includegraphics[width=0.8 \linewidth]{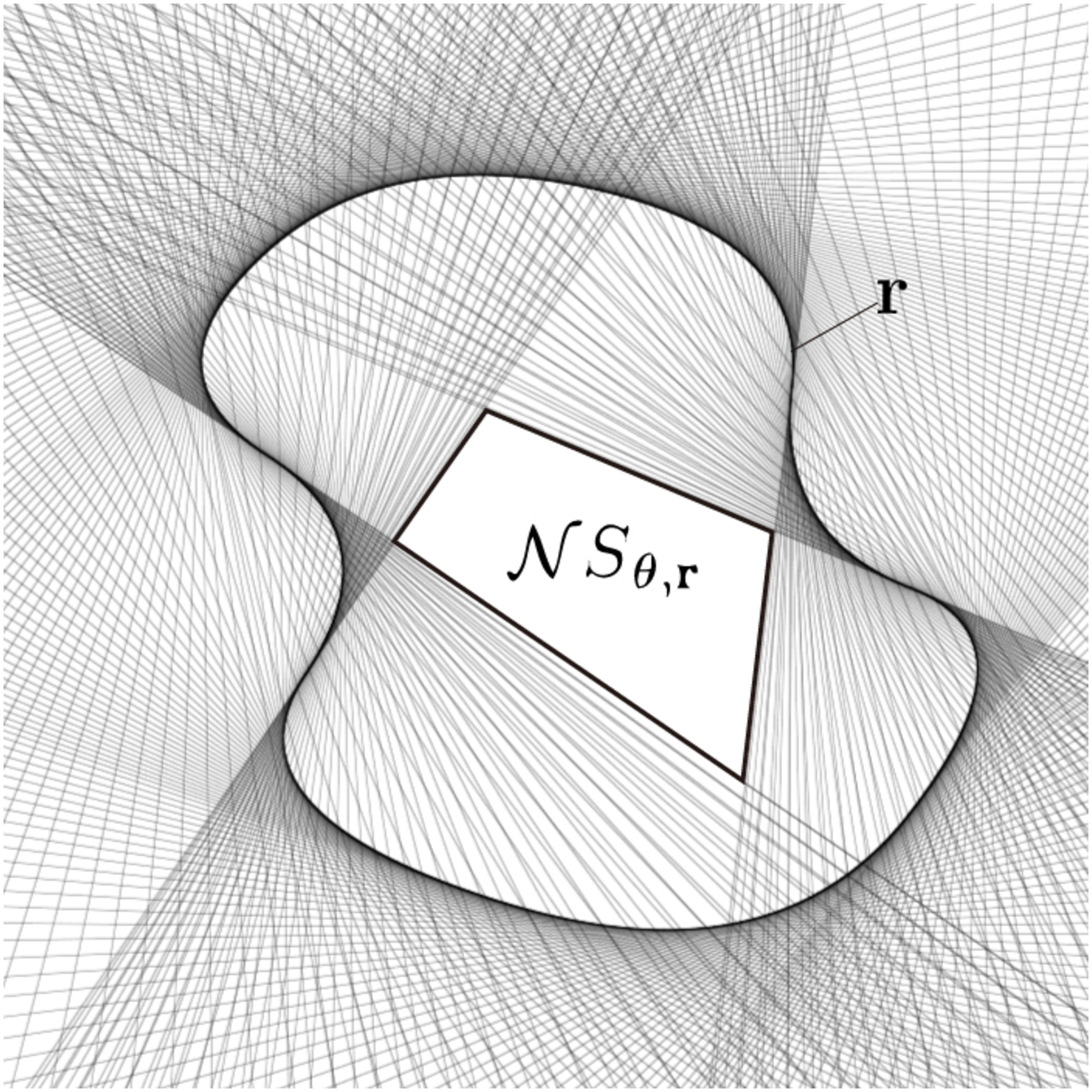} 
	\end{center}
  \end{minipage}
  \hspace{2.0pc} 
  \begin{minipage}{.45\linewidth} 
  	\begin{center}
\includegraphics[width=0.8 \linewidth]{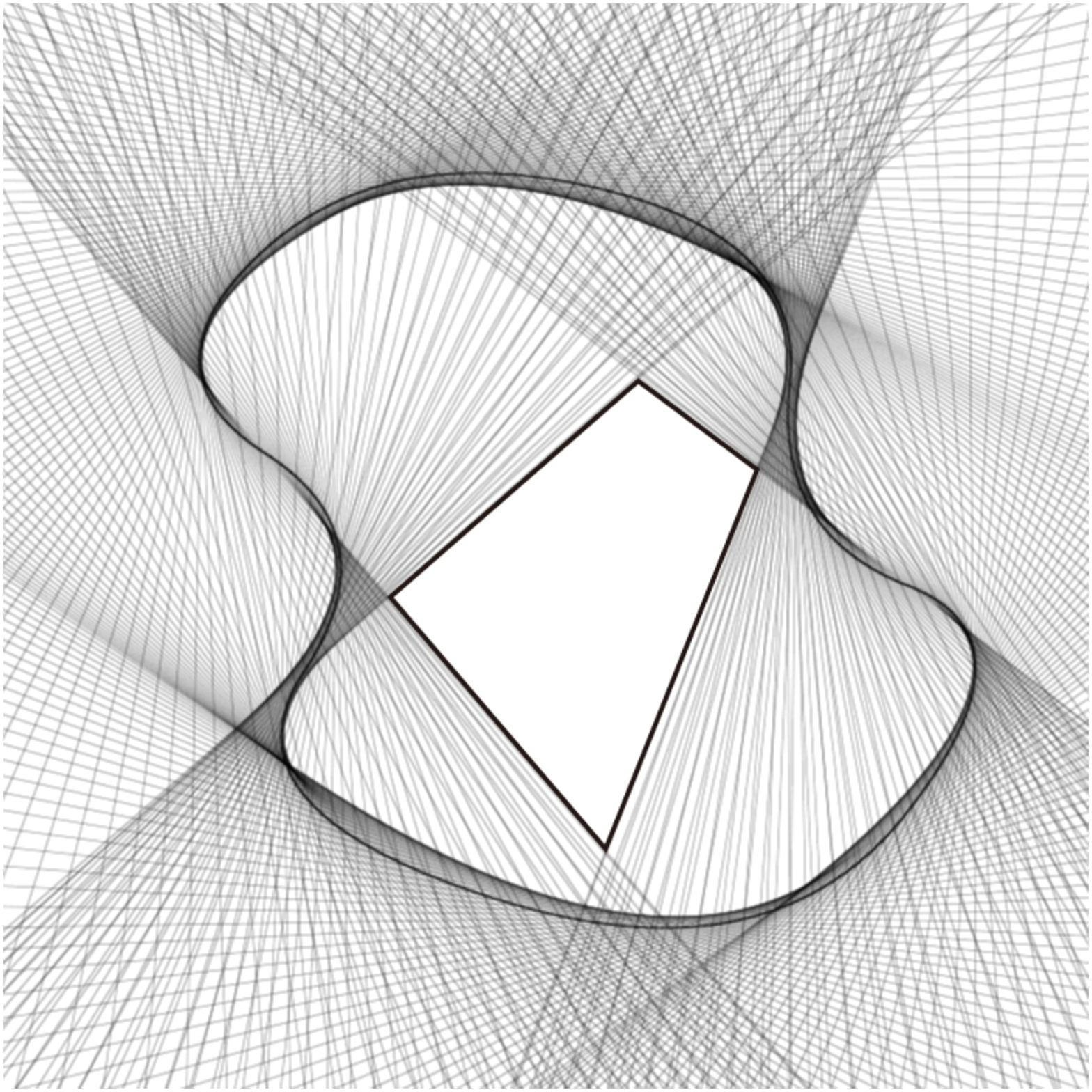} 
	\end{center}
  \end{minipage}
  \end{center}

  \begin{center}
  \begin{minipage}{.45\linewidth} 
  \begin{center}
  \includegraphics[width=0.8 \linewidth]{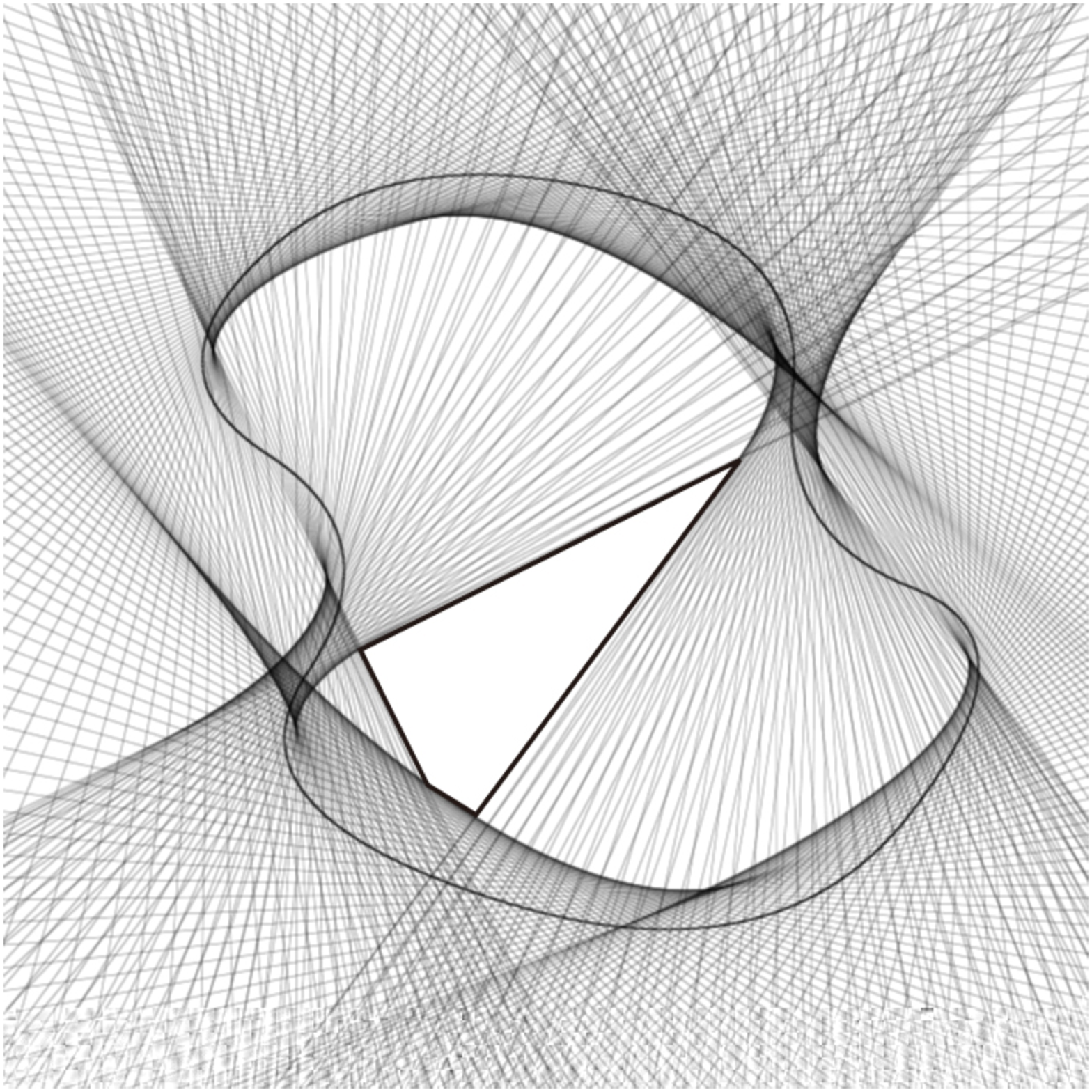} 
	\end{center}
  \end{minipage}
  \hspace{2.0pc} 
  \begin{minipage}{.45\linewidth} 
  \begin{center}
 \includegraphics[width=0.8 \linewidth]{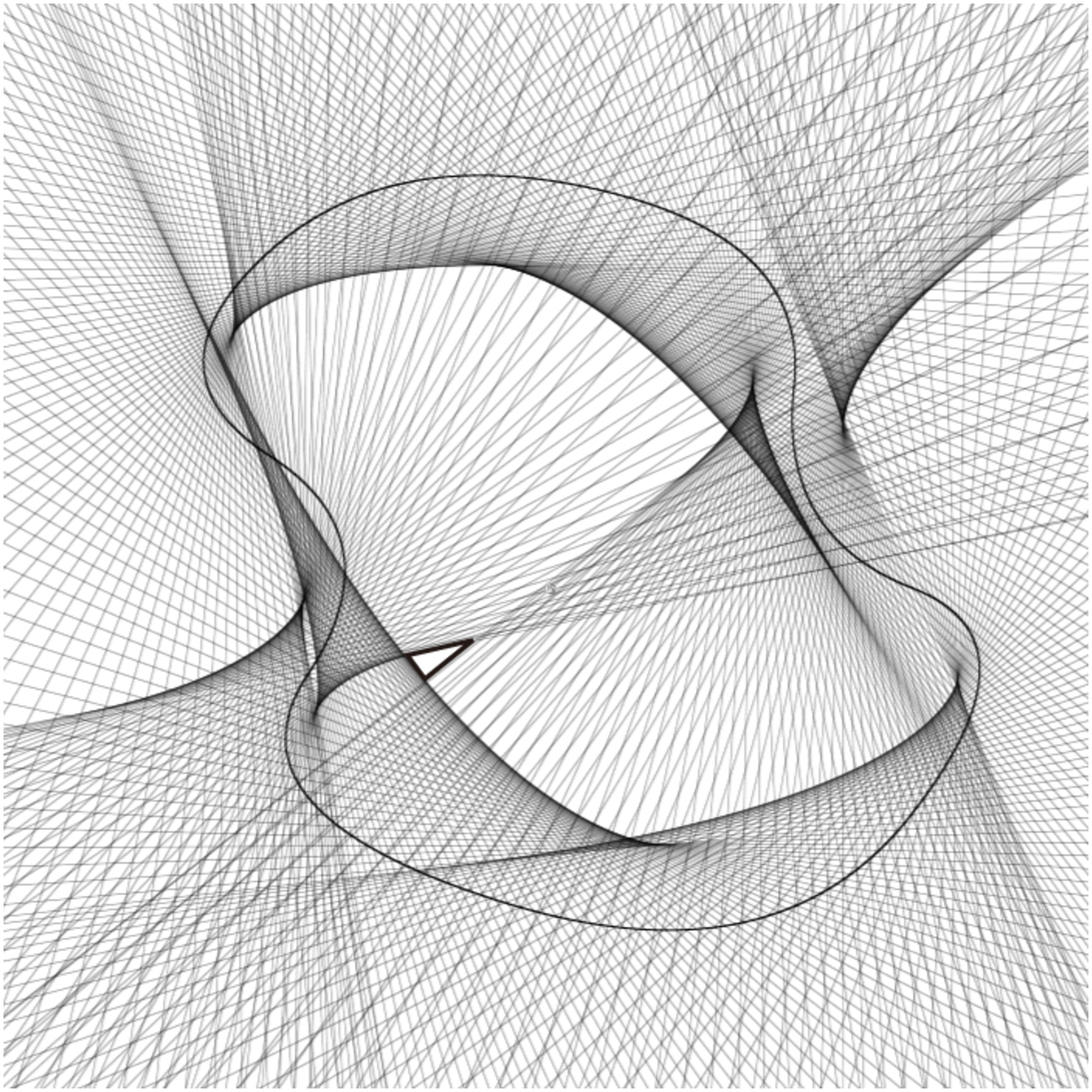} 
	\end{center}
  \end{minipage}
  \end{center}
  \caption{${\mathcal NS}_{\theta, {\bf r}}$ for several $\theta$s.    
Left top : $\theta = 0$, right top : $\theta = \pi/12$, left bottom :  $\theta = \pi/6$, right bottom : $\theta = \pi/4$.}
  \label{figure 2}
  \end{figure}
When the given ${\bf r}$ has its no-silhouette $\mathcal{NS}_{\bf r}$, 
by definition, it follows that $\mathcal{NS}_{{\bf r}}=\mathcal{NS}_{0, {\bf r}}$.  
\begin{lemma}\label{lemma 1.1}
For any $C^\infty$ immersion 
${\bf r}: S^1\to \mathbb{R}^2$, 
$\mathcal{NS}_{\frac{\pi}{2}, {\bf r}}$ is the empty set.  
\end{lemma}
\underline{\it Proof of Lemma \ref{lemma 1.1}}\qquad 
For any point $P \in {\mathbb R}^2$, let $F_P: S^1\to \mathbb{R}$ be the 
function defined by  
$$ 
F_P(s) = (P - {\bf r}(s)) \cdot (P - {\bf r}(s)), 
\leqno{(1.2)}
\label{(1.2)}
$$
where the dot in the center stands for the scalar product of two vectors.       
Since $F_P$ is a $C^\infty$ function and $S^1$ is compact, 
there exist the maximum and the minimum of the set of images 
$\left\{F_P (s) \mid s \in S^1 \right\}$.
Let $s_1$ (resp., $s_2$) be a point of $S^1$ 
at which $F_P$ attains its maximum (resp., minimum).   
Then, both $s_1$ and $s_2$ are critical points of $F_P$.    
Thus, differentiating (1.2) with respect to $s$ yields that 
the vector $(P - {\bf r}(s_i))$ is perpendicular to the tangent line to ${\bf r}$ at ${\bf r}(s_i)$. 
It follows that 
$P \in \left({\bf r}(s_i)+R_\frac{\pi}{2} (d{\bf r}_{s_i}(T_{s_i}S^1)\right)$.
\hfill\qquad $\Box$ 
\par 
\medskip    
In Section \ref{section 2}, 
it turns out that with respect to the Pompeiu-Hausdorff metric 
the topological closure of $\mathcal{NS}_{\theta, {\bf r}}$ 
varies continuously depending on 
$\theta$ while $\mathcal{NS}_{\theta, {\bf r}}$ is not empty (Proposition \ref{proposition 5}).       
Therefore, by Lemma 1.1, the following notion of aperture angle $\theta_{\bf r}$ $(0<\theta_{\bf r}\le \frac{\pi}{2})$ 
is well-defined.  
\begin{definition}\label{definition 1}
{\rm 
Let ${\bf r}: S^1\to \mathbb{R}^2$ be a 
$C^\infty$ immersion with its no-silhouette $\mathcal{NS}_{\bf r}$.   
Then, $\theta_{\bf r}$ $(0 < \theta_{\bf r}\le \frac{\pi}{2})$ 
is defined as the largest angle which satisfies $\mathcal{NS}_{\theta, {\bf r}}\ne \emptyset$ 
for any $\theta$ ($0\le \theta< \theta_{\bf r}$).    
The angle $\theta_{\bf r}$ is called the {\it aperture angle} of the given ${\bf r}$.   
}
\end{definition} 
In Section 2, it turns out also that $\overline{\mathcal{NS}_{\theta, {\bf r}}}$ is a Wulff shape 
for any $\theta$ $(0\le \theta <\theta_{\bf r})$, where 
$\overline{\mathcal{NS}_{\theta, {\bf r}}}$ stands for the topological closure of 
$\mathcal{NS}_{\theta, {\bf r}}$ (Proposition \ref{proposition 4}).      
We are interested in how the Wulff shape $\overline{\mathcal{NS}_{\theta, {\bf r}}}$ 
dissolves   
as $\theta$ goes to $\theta_{\bf r}$ 
from $0$.       
\begin{theorem}\label{theorem 1}
Let ${\bf r}: S^1\to \mathbb{R}^2$ be a $C^\infty$ immersion with its no-silhouette 
$\mathcal{NS}_{\bf r}$.    Then, for any $\theta$ $(0<\theta<\theta_{\bf r})$, 
$\overline{\mathcal{NS}_{\theta}, {\bf r}}$ is never a polygon even if 
the given 
$\overline{\mathcal{NS}_{\bf r}}$ is a polygon.    
\end{theorem}
\noindent 
By Theorem \ref{theorem 1}, none of $\overline{\mathcal{NS}_{\frac{\pi}{12}, {\bf r}}}$, 
$\overline{\mathcal{NS}_{\frac{\pi}{6}, {\bf r}}}$ and $\overline{\mathcal{NS}_{\frac{\pi}{4}, {\bf r}}}$ in 
Figure \ref{figure 2} is a polygon although $\overline{\mathcal{NS}_{0, {\bf r}}}$ is a polygon constructed 
by four tangent lines to ${\bf r}$ at four inflection points.       
\begin{theorem}\label{theorem 2}
\textcolor{black}{
Let ${\bf r}: S^1\to \mathbb{R}^2$ be a $C^\infty$ immersion with its no-silhouette 
$\mathcal{NS}_{\bf r}$.    
Then, there exists the unique point $P_{\bf r}\in \mathbb{R}^2$ such that 
for any sequence $\{\theta_i\}_{i=1, 2, \ldots}\subset [0, \theta_{\bf r})$ 
satisfying $\lim_{i\to \infty}\theta_i =\theta_{\bf r}$, the following holds:   .    
\[
\lim_{i\to \infty}d_H(\overline{\mathcal{NS}_{\theta_i}, {\bf r}}, P_{\bf r})=0. 
\]   .  
}
\end{theorem}
\noindent 
Here, $d_H: \mathcal{H}(\mathbb{R}^2)\times \mathcal{H}(\mathbb{R}^2)\to \mathbb{R}$ is the
Pompeiu-Hausdorff metric (for the Pompeiu-Hausdorff metric, 
see Section 2).  
Theorem \ref{theorem 2} justifies the following definition.  
\begin{definition}\label{definition 2}
{\rm 
Let ${\bf r}: S^1\to \mathbb{R}^2$ be a $C^\infty$ immersion with its no-silhouette 
$\mathcal{NS}_{\bf r}$.   Then, 
the set $\cup_{\theta\in[0,\theta_{\bf r})}\overline{\mathcal{NS}_{\theta, {\bf r}}}$ 
is called the {\it aperture} of ${\bf r}$ and the unique point 
$P_{\bf r}=\lim_{\theta\to \theta_{\bf r}}\overline{\mathcal{NS}_{\theta, {\bf r}}}$ is called the {\it aperture point} of {\bf r}.  
Here, $\theta_{\bf r}$ $(0 < \theta_{\bf r}\le \frac{\pi}{2})$ is the aperture angle of ${\bf r}$.  
}  
\end{definition}
\noindent 
The simplest example is a circle.   
The aperture of a circle is the topological closure of its inside region 
and the aperture point of it is its center.   
In this case, the aperture angle is $\pi/2$.         
In general, in the case of curves with no inflection points,  
the crystal growth is relatively simpler than \textcolor{black}{in the case of} curves with inflections as follows.   
\begin{theorem}\label{theorem 3}
Let ${\bf r}: S^1\to \mathbb{R}^2$ be a $C^\infty$ immersion with its no-silhouette 
$\mathcal{NS}_{\bf r}$.   
Suppose that ${\bf r}$ has no inflection points.   
Then, for any two $\theta_1, \theta_2$ satisfying $0\le \theta_1<\theta_2<\theta_{\bf r}$, 
the following inclusion holds:    
\[
\mathcal{NS}_{\theta_1, {\bf r}}\supset 
\mathcal{NS}_{\theta_2, {\bf r}}.   
\]
\end{theorem}
\noindent 
Figure \ref{figure 2} shows that in general it is impossible to expect the same property for a curve with 
inflection points. 
\bigskip 
\par 
In Section 2, preliminaries are given.         
Theorems \ref{theorem 1}, \ref{theorem 2} and \ref{theorem 3} 
are proved in Sections \ref{section 3}, \ref{section 4} and \ref{section 5} 
respectively. 
 \section{Preliminaries} \label{section 2}
\subsection{Spherical curves}
Let $\widetilde{\bf r}: S^1\to S^2$ be a $C^\infty$ immersion.     
Let $\widetilde{\bf t}: S^1\to S^2$ be the mapping defined by 
\[
\widetilde{\bf t}(s)=
\frac{\widetilde{\bf r}'(s)}{||\widetilde{\bf r}'(s)||}, 
\]
where $\widetilde{\bf r}'(s)$ stands for 
differentiating $\widetilde{\bf r}(s)$ 
with respect to $s\in S^1$.    
Let $\widetilde{\bf n}: S^1\to S^2$ be the mapping defined by  
\[
\det \left(\widetilde{\bf r}(s), \widetilde{\bf t}(s), 
\widetilde{\bf n}(s)\right)=1.   
\]
The mapping $\widetilde{\bf n}: S^1\to S^2$ is called the {\it spherical dual of}\/ $\widetilde{\bf r}$.            
The singularities of $\widetilde{\bf n}$ belong to the class of Legendrian 
singularities which are relatively well-investigated (for instance, see \cite{arnold,arnolddynamical8, arnoldetall}).    
Let $U$ be an open arc of $S^1$.   
Suppose that $||\widetilde{\bf r}'(s)||=1$ for any $\textcolor{black}{s}\in U$.    
Then, for the orthogonal moving frame $\{{\bf r}(s), {\bf t}(s), 
{\bf n}(s)\}$, $(s\in U)$,  
the following Serre-Frenet type formula has been known.   
\begin{lemma}[\cite{pedal, demonstratio}]\label{lemma 2.1}
\[
\left\{
\begin{array}{ccl}
\widetilde{\bf r}'(s) & = & \widetilde{\bf t}(s) \\ 
\widetilde{\bf t}'(s) & = & -\widetilde{\bf r}(s)
+\kappa_g(\theta)\widetilde{\bf n}(s) \\ 
\widetilde{\bf n}'(s) & = & -\kappa_g(\theta)\widetilde{\bf t}(s).  
\end{array}
\right.
\]
\end{lemma}
\noindent 
Here, $\kappa_g(\theta)$ is defined by 
\[
\kappa_g(\theta)=\det \left(\widetilde{\bf r}(s), \widetilde{\bf t}(s), 
\widetilde{\bf t}'(s)\right).
\]
\par 
\medskip 
Let $N$ be the north pole $(0,0,1)$ of the unit sphere $S^2\subset \mathbb{R}^3$ and 
let $S^2_{N,+}$ be the northern hemisphere $\{P\in S^2\; |\; N\cdot P>0\}$, 
where $N\cdot P$ stands for the scalar product of two vectors $N, P\in \mathbb{R}^3$.  
Then, define the mapping $\alpha_N: S^2_{N,+}\to \mathbb{R}^2\times \{1\}$, which is called 
the {\it central projection}, as follows:  
\[
\alpha_N(P_1, P_2, P_3)=\left(\frac{P_1}{P_3}, \frac{P_2}{P_3}, 1\right),  
\]
where $P=(P_1, P_2, P_3)\in S^2_{N,+}$.         
Let ${\bf r}: S^1\to \mathbb{R}^2$ be a $C^\infty$ immersion.   
Then, from ${\bf r}$ 
we can naturally obtain a spherical curve $\widetilde{\bf r}: S^1\to S^2$ as follows:
\[
\widetilde{\bf r}=\alpha_N^{-1}\circ Id \circ {\bf r},
\] 
where $Id: \mathbb{R}^2\to \mathbb{R}^2\times \{1\}$ is the mapping defined by 
$Id(P)=(P, 1)$.   
For any $s\in S^1$, let $GC_{\widetilde{\bf r}(s)}$ be the intersection 
$(\mathbb{R}\widetilde{\bf r}(s)+\mathbb{R}\widetilde{\bf t}(s))\cap S^2$.    
The following clearly holds:  
\begin{lemma}\label{lemma 2.2}
By the central projection $\alpha_N: S^2_{N,+}\to \mathbb{R}^2\times \{1\}$, 
$GC_{\widetilde{\bf r}(s)}\cap S^2_{N,+}$ is mapped to the line ${\bf r}(s)+d{\bf r}_s(T_s(S^1))$.   
\end{lemma}
One of the merit of considering inside the sphere $S^2$ is the following:   
\begin{lemma}[\cite{nishimurasakemi1}]\label{lemma 2.3}
Let $\textcolor{black}{\widetilde{\bf r}}: S^1\to S^2$ be a Legendrian mapping.    
Then, the following two are equivalent conditions.   
\begin{enumerate}
\item The set 
$$
S^2-\bigcup_{s\in S^2}GC_{\widetilde{r}(s)}
$$
is not empty and $N$ is inside this open set.   
\item The connected subset 
$\{\widetilde{\bf n}(s)\; |\; s\in S^1\}$ is inside $S^2_{N,+}$, 
where $\widetilde{\bf n}$ is the dual of $\widetilde{\bf r}$. 
\end{enumerate}
\end{lemma}
Let $\Psi_N:S^2-\{\pm N\}\to S^2$ be the mapping defined by 
$$
\Psi_N(P)=\frac{1}{\sqrt{1-(N\cdot P)^2}}(N-(N\cdot P)P).    
$$
The mapping $\Psi_N$ is very useful for studying spherical pedals, pedal unfoldings of spherical pedals, 
hedgehogs, and Wulff shapes (see \cite{pedal, demonstratio, nishimuraunfolding, nishimurasakemi1, nishimurasakemi2}). 
There is also a hyperbolic version of $\Psi_N$ (\cite{izumiyatari}).    
The fundamental properties of $\Psi_N$ is as follows:    
\begin{enumerate}
\item For any $P\in S^2-\{\pm N\}$, the equality $P\cdot \Psi_N(P)=0$ holds,    
\item for any $P\in S^2-\{\pm N\}$, the property $\Psi_N(P)\in \mathbb{R}N+\mathbb{R}P$ holds,     
\item for any $P\in S^2-\{\pm N\}$, the property $N\cdot \Psi_N(P)>0$ holds,      
\item the restriction $\Psi_N|_{S^2_{N,+}-\{N\}}: S^2_{N,+}-\{N\}\to S^2_{N,+}-\{N\}$ is a $C^\infty$ diffeomorphism.   
\end{enumerate}
By these properties, we have the following: 
\begin{lemma}
The mapping 
$\alpha_N\circ \Psi_N\circ \alpha_N^{-1} \textcolor{black}{: \mathbb{R}^2\times \{1\}-\{N\}\to \mathbb{R}^2\times \{1\}-\{N\}}$ 
is the inversion of $\mathbb{R}^2\times \{1\}- \{N\}$ \textcolor{black}{with respect to $N$}.       
\end{lemma}
\subsection{Spherical polar sets and the spherical polar transform}
\begin{figure}[hbtp]
\begin{center}
\includegraphics[width=4cm]{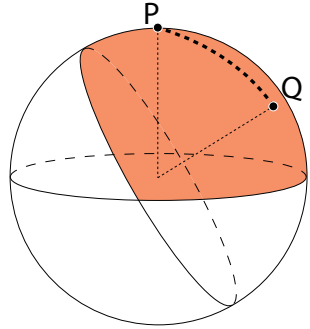}
\caption{Spherical polar set $\{P, Q\}^\circ=(PQ)^\circ$.}
\label{figure 3}
\end{center}
\end{figure}

\begin{figure}[hbtp]
\begin{center}
\includegraphics[width=4cm]{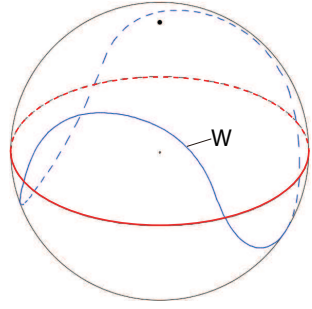}
\caption{Not hemispherical $W\subset S^2$.}
\label{figure 4}
\end{center}
\end{figure}

For any point $P$ of $S^2$, we let $H(P)$ be the following set:
$$
H(P)=\{Q\in S^2\; |\; P\cdot Q\ge 0\}.
$$
Here, the dot in the center stands for the scalar product of $P, Q\in \mathbb{R}^3$.    
\begin{definition}[\cite{nishimurasakemi2}]\label{spherical polar set}
{\rm 
Let $W$ be a subset of $S^2$.   Then, the set 
$$\bigcap_{P\in W}H(P)$$ 
is called the {\it spherical polar set} of $W$ and is denoted by $W^\circ$.   
}
\end{definition}  
Figure \ref{figure 3} illustrates Definition \ref{spherical polar set}.   
It is clear that $W^\circ=\cap_{P\in W}H(P)$ is closed for any $W\subset S^2$.      
\begin{definition}[\cite{nishimurasakemi2}]\label{definition 2.1}
{\rm 
A subset $W\subset S^2$ is said to be {\it hemispherical} if there exists a point $P\in S^2$ such that 
$H(P)\cap W=\emptyset$.   
}
\end{definition}
Figure \ref{figure 4} illustrates Definition \ref{definition 2.1}.  
\begin{definition}[\cite{nishimurasakemi2}]\label{definition 2.2}
{\rm 
A hemispherical subset $W\subset S^2$ is said to be {\it spherical convex} if $PQ\subset W$ for any $P, Q\in W$.     
}
\end{definition}
Here, $PQ$ stands for the following arc:  
$$
PQ=\left\{\left.\frac{(1-t)P+tQ}{||(1-t)P+tQ||}\in S^2\; \right|\; 0\le t\le 1\right\}.  
$$
Note that $||(1-t)P+tQ||\ne 0$ for any $P, Q\in W$ and any $t\in [0,1]$ if $W\subset S^2$ is 
hemispherical.       
Note further that $W^\circ$ is spherical convex if $W$ is hemispherical and it has an interior point.   
\begin{definition}[\cite{nishimurasakemi2}]\label{definition 2.2}
{\rm Let $W$ be a hemispherical subset of $S^2$.     
Then, the {\it spherical convex hull of} $W$ (denoted by 
$\mbox{\rm s-conv}(W)$) is the following set.   
$$
\mbox{\rm s-conv}(W)= 
\left\{\left.
\frac{\sum_{i=1}^k t_iP_i}{||\sum_{i=1}^kt_iP_i||}\;\right|\; 
P_i\in W,\; \sum_{i=1}^kt_i=1,\; t_i\ge 0, k\in \mathbb{N}
\right\}.
$$ 
}
\end{definition}
\begin{lemma}[Lemma 2.5 of \cite{nishimurasakemi2}]\label{Maehara's lemma}
For any hemispherical finite subset $W=\{P_1, \ldots, P_k\}\subset S^{n+1}$, the following holds:
$$
\left\{\left.
\frac{\sum_{i=1}^k t_iP_i}{||\sum_{i=1}^kt_iP_i||}\;\right|\; 
P_i\in W,\; \sum_{i=1}^kt_i=1,\; t_i\ge 0
\right\}^\circ 
= 
H(P_1)\cap \cdots \cap H(P_k).
$$
\end{lemma}
\noindent 
Lemma \ref{Maehara's lemma} is called {\it Maehara's lemma} (see \cite{nishimurasakemi2}).   
\par 
\medskip 
\begin{definition}[\cite{barnsley}]\label{Pompeiu-Hausdorff}
{\rm Let $(X, d)$ be a complete metric space.     
\begin{enumerate}
\item Let $x$ be a point of $X$ and let $B$ a non-empty compact subset 
of $X$.     Define 
\[
d(x,B)=\min\{d(x,y)\; |\; y\in B\}.
\] 
Then, $d(x,B)$ is called the {\it distance from the point $x$ to the set $B$}.   
\item Let $A, B$ be two non-empty compact subsets of $X$.   
Define 
\[
d(A,B)=\max \{d(x, B)\; |\; x\in A\}.   
\] 
Then, $d(A,B)$ is called the {\it distance from the set $A$ to the set $B$}.  
\item Let $A, B$ be two non-empty compact subsets of $X$.   
Define 
\[
d_H(A,B)=\max \{d(A, B), d(B, A)\}.   
\] 
Then, $d_H(A,B)$ is called the {\it Pompeiu-Hausdorff distance between $A$ and $B$}.  
\end{enumerate}
}
\end{definition}
\par 
\noindent 
Let $(X, d)$ be a complete metric space.     
The set consisting of non-empty compact subsets of $X$ is denoted by 
$\mathcal{H}(X)$, which is the metric space with respect to 
the Pompeiu-Hausdorff metric $d_H: \mathcal{H}(X)\times\mathcal{H}
(X)\to \mathbb{R}_+\cup\{0\}$, where $d_H$ is 
the metric naturally induced by the Pompeiu-Hausdorff distance.  
It is well-known also that the metric space $(\mathcal{H}(X), d_H)$ is complete.     
For more details on the complete metric space 
$(\mathcal{H}(X), d_H)$, see for instance \cite{barnsley, falconer}.      
\begin{definition}\label{spherical polar transform}
{\rm 
Let $\bigcirc: \mathcal{H}(S^2)\to \mathcal{H}(S^2)$ be the mapping defined by 
\[
\bigcirc(A)=A^\circ.   
\]
}
\end{definition}
\par 
\noindent 
The mapping $\bigcirc: \mathcal{H}(S^2)\to \mathcal{H}(S^2)$ is called 
the {\it spherical polar transform}.   
\begin{proposition}\label{continuity}
The spherical polar transform 
is continuous with respect to 
the Pompeiu-Hausdorff metric.    
\end{proposition}
\underline{\it Proof of Proposition \ref{continuity}}\qquad 
Let $\{A_i\}_{i=1, 2, \ldots}\subset \mathcal{H}(S^2)$ be a convergent sequence.   
Set $A=\lim_{i\to \infty} A_i$.     In order to prove Proposition \ref{continuity}, 
it is sufficient to show that $A^\circ=\lim_{i\to \infty}A_i^\circ$.    
\par 
Suppose that there exists a real number $\varepsilon >0$ such that 
for any $n\in \mathbb{N}$ there exists an $i_n$ $(i_n > n)$ such that 
$d_H(A^\circ_{i_n}, A^\circ) > \varepsilon$.    Then, by Definition \ref{Pompeiu-Hausdorff}, 
it follows that for any $n\in \mathbb{N}$, at least one of 
$d(A^\circ_{i_n}, A^\circ) > \varepsilon$ and $d(A^\circ, A^\circ_{i_n}) > \varepsilon$ holds.    By taking a subsequence if necessary, from the first 
we may assume that one of the following holds:   
\begin{enumerate}
\item $d(A^\circ_{i_n}, A^\circ) > \varepsilon$ for any $n\in \mathbb{N}$.  
\item $d(A^\circ, A^\circ_{i_n}) > \varepsilon$ for any $n\in \mathbb{N}$, 
\end{enumerate}   
\par 
\smallskip 
We first show that (1) implies a contradiction.   
By Definition \ref{Pompeiu-Hausdorff}, it follows that for any $n\in \mathbb{N}$ there exists a point 
$x_n\in A^\circ_{i_n}$ such that $d(x_n, A^\circ) > \varepsilon$.    
Again by Definition \ref{Pompeiu-Hausdorff}, it follows that for any $n\in \mathbb{N}$ 
there exists a point  $x_n\in A^\circ_{i_n}$ such that  the inequality 
$d(x_n, y)>\varepsilon$ holds for any $y\in A^\circ$.     It is known that $A$ can be characterized as follows 
(\cite{barnsley}).   
\[
A=\left\{
P\in S^2\; |\; \exists P_n\in A_{i_n} (n\in \mathbb{N})\mbox{ such that }
\lim_{n\to \infty}P_n=P
\right\}.
\leqno{(2.1)}
\]
Let $P$ be a point of $A$.    
By (2.1), for any $n\in \mathbb{N}$ we may choose a point $P_n\in A_{i_n}$ such that 
$\lim_{n\to \infty}P_n=P$.    Then, since $x_n\in A_{i_n}^\circ$, it follows that 
$x_n\cdot P_n\ge 0$.     
Since $S^2$ is compact, there exists a convergent 
subsequence $\{x_{j_n}\}_{n=1, 2, \ldots}$ of the sequence $\{x_n\}_{n=1, 2, \ldots}$.   
Set $x=\lim_{n\to \infty}x_{j_n}$.   Then, the inequality 
$d(x_n, y)>\varepsilon$ implies the inequality $d(x, y)\textcolor{black}{\ge}\varepsilon$ for any 
$y\in A^\circ$.    On the other hand, the inequality $x_n\cdot P_n\ge 0$ 
implies the inequality $x\cdot P\ge 0$ for any 
$P\in A$.    Therefore, the point $x$ is an element of $A^\circ$ 
such that the inequality  
$d(x, y)\textcolor{black}{\ge}\varepsilon$ holds for any $y\in A^\circ$.   This is a contradiction.   
\par 
\smallskip 
We next show that (2) implies a contradiction.    
By the same argument as in (1), we have that 
for any $n\in \mathbb{N}$ 
there exists a point  $x_n\in A^\circ$ such that  the inequality 
$d(x_n, y_n)>\varepsilon$ for any $y_n\in A^\circ_{i_n}$.    
This implies that there exists an $M\in \mathbb{N}$ such that 
for any $n\in \mathbb{N}$ there exists $P_n\in \textcolor{black}{A_{i_n}}$ such that 
$x_n\cdot P_n < -\frac{\varepsilon}{M}$.    Since $S^2$ is compact, there exists 
a subsequence $\{j_n\}_{n=1, 2, \ldots}$  of $\mathbb{N}$ such that both 
$\{x_{j_n}\}_{n=1, 2, \ldots}$, $\{P_{j_n}\}_{n=1, 2, \ldots}$ are convergent sequences.    
Set $x=\lim_{n\to \infty}x_{j_n}$ and $P=\lim_{n\to \infty}P_{j_n}$.      
Then, the inequality $x_n\cdot P_n < -\frac{\varepsilon}{M}$ implies 
the inequality $x\cdot P \le -\frac{\varepsilon}{M}$.   
On the other hand, since $A^\circ$ is compact, $x$ belongs to $A^\circ$.    
Moreover, by (2.1), $P$ belongs to $A$.      Hence, by Definition \ref{spherical polar set}, 
the scalar product 
$x\cdot P$ must be non-negative.   
Therefore, we have a contradiction.   
\hfill\qquad $\Box$ 
\subsection{Wulff shapes}
\begin{figure}[hbtp]
\begin{center}
\includegraphics[width=6cm]{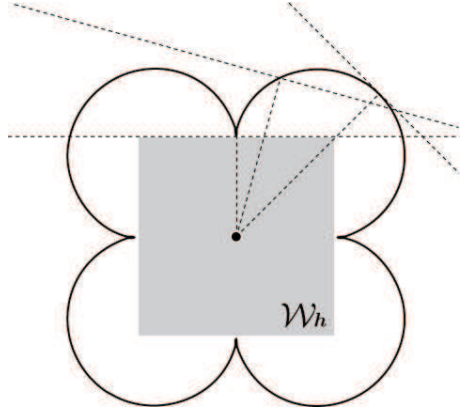}
\caption{The Wulff shape associated with the support function $h$.}
\label{figure 5}
\end{center}
\end{figure}

Let $\mathbb{R}_+$ be the set $\{\lambda\in \mathbb{R}\; |\; \lambda>0\}$ 
and let $h:S^1\to \mathbb{R}_+$ be a continuous function.  
For any $s \in S^1\subset \mathbb{R}^2$, set  
$$
\Gamma_{h,s}
=\{P\in \mathbb{R}^{2}\; |\; P\cdot s\le h(s)\},  
$$
where the dot in the center stands for the scalar product of two vectors $P, s\in \mathbb{R}^2$.  
The following set is called the 
{\it Wulff shape associated with the support function $h$} (see Figure \ref{figure 5}):   
$$
\mathcal{W}_h=\bigcap_{s\in S^1}\Gamma_{h, s}.
$$
For any crystal at equilibrium the shape of it can be constructed as 
the Wulff shape $\mathcal{W}_{h}$ 
by an appropriate support function $h$ (\cite{wulff}).        
It is clear that any Wulff shape $\mathcal{W}_h$ is a convex body (namely, it is compact, convex and the origin of 
$\mathbb{R}^{\textcolor{black}{2}}$ is contained in $\mathcal{W}_h$ as an interior point).       
It has been known that its converse, too, holds as follows.    
\begin{proposition}[p.~573 of  \cite{taylor}]\label{taylor}
Let $W$ be a subset of $\mathbb{R}^{2}$.    
Then, there exists a parallel translation 
$T: \mathbb{R}^{2}\to \mathbb{R}^{2}$ such that 
$T(W)$ is the Wulff shape associated with an appropriate support function if and only if 
$W$ is a convex body.
\label{proposition 1} 
\end{proposition}
\begin{proposition}[Theorem 1.1 of \cite{nishimurasakemi2}]\label{proposition 2}
Let $\left\{\mathcal{W}_{h_i}\right\}_{i=1, 2, \ldots}$ 
be a Cauchy sequence of Wulff shapes in 
$\mathcal{H}_{\mbox{\rm conv}}(\mathbb{R}^{2})$ 
with respect to the Pompeiu-Hausdorff metric $d_H$.    
Suppose that $\lim_{i\to \infty}\mathcal{W}_{h_i}$ 
does not have an interior point.     
Then, it must be a point or a segment.  
\end{proposition} 
\begin{proposition}[Theorem 1.2 of \cite{nishimurasakemi2}]\label{dual Wulff shape}
Let $h: S^1\to \mathbb{R}_+$ be a continuous function.     
Then, for the Wulff shape $\mathcal{W}_h$, the set 
$Id^{-1}\circ\alpha_N\left(\left(\alpha_N^{-1}\circ Id(\mathcal{W}_h)\right)^\circ\right)$  
is the Wulff shape associated with an appropriate support function.   
\end{proposition}
The Wulff shape  
$Id^{-1}\circ\alpha_N\left(\left(\alpha_N^{-1}\circ Id(\mathcal{W}_h)\right)^\circ\right)$ 
is called the {\it dual Wulff shape of} $\mathcal{W}_h$.   
\begin{proposition}[Theorem 1.3 of \cite{nishimurasakemi2}]\label{proposition C1}
Let $h: S^1\to \mathbb{R}_+$ be a function of class $C^1$.     
Then, the Wulff shape $\mathcal{W}_h$ is never a polygon.   
\label{proposition 3}
\end{proposition}  
\begin{proposition}\label{proposition 4}
Let ${\bf r}: S^1\to \mathbb{R}^2$ be a $C^\infty$ immersion with its no-silhouette 
$\mathcal{NS}_{\bf r}$.    
Then, for any $\theta\in [0,\theta_{\bf r})$, 
there exists a parallel translation $T_\theta: \mathbb{R}^2\to \mathbb{R}^2$ such that 
$T_\theta(\overline{\mathcal{NS}_{\theta, {\bf r}}})$ is a Wulff shape $\mathcal{W}_{h_\theta}$ 
by an appropriate support function 
$h_\theta: S^1\to \mathbb{R}_+$.    
\end{proposition}
\underline{\it Proof of Proposition \ref{proposition 4}}\qquad 
We first show that $\mathcal{NS}_{\theta, {\bf r}}$ is an open set for any $\theta\in [0,\theta_{\bf r})$.      
Let $P$ be a point of $\mathcal{NS}_{\theta, {\bf r}}$.    
Suppose that for any positive integer $n$, 
there exists a point $P_n\in D(P, \frac{1}{n})\cap\left(\cup_{s\in S^1}\left({\bf r}(s)+
R_\theta \left(d{\bf r}_s(T_s(S^1))\right)\right)\right)$, 
where $D(P, \frac{1}{n})$ is the disc $D(P, \frac{1}{n})=\{Q\in \mathbb{R}^2\; |\; ||P-Q||\le\frac{1}{n}\}$.      
Then, since $S^1$ is compact, by taking a subsequence if necessary, we may assume that 
there exists a convergent sequence $s_n\in S^1$ $(n\in \mathbb{N})$ such that $P_n$ belongs to 
$D(P, \frac{1}{n})\cap\left({\bf r}(s_n)+R_\theta \left(d{\bf r}_{s_n}(T_{s_n}(S^1))\right)\right)$.   
Then, we have that $P\in {\bf r}(s)+R_\theta \left(d{\bf r}_{s}(T_{s}(S^1))\right)$ where $s=\lim_{i\to \infty}s_i$, which implies 
$P\not\in \mathcal{NS}_{\theta, {\bf r}}$.      Hence,  $\mathcal{NS}_{\theta, {\bf r}}$ is an open set.   
\par 
Since $\theta<\theta_{\bf r}$, 
it follows that $\mathcal{NS}_{\theta, {\bf r}}\ne \emptyset$.   
Let $P$ be a point of $\mathcal{NS}_{\theta, {\bf r}}$.   
Let $P_s\in {\bf r}(s)+R_\theta d{\bf r}_s(T_s(S^1))$ 
be the point such that the vector $PP_s$ is perpendicular to the line 
${\bf r}(s)+R_\theta d{\bf r}_s(T_s(S^1))$.  
Then, by obtaining the concrete expression of $P_s$, 
it follows that the mapping $f: S^1\to \mathbb{R}^2$ defined by $f(s)=P_s$ is of class $C^\infty$.        
\textcolor{black}
{
By Subsection 2.1 and \cite{pedal}, 
the mapping $f: S^1\to \mathbb{R}^2$ is exactly the pedal curve 
of the family of lines $\left\{
{\bf r}(s)+R_\theta d{\bf r}_s(T_s(S^1))
\right\}_{s\in S^1}$ relative to the pedal point $P\in \mathcal{NS}_{\theta, {\bf r}}$. 
}
Let $I: \mathbb{R}^2-\{P\}\to \mathbb{R}^2-\{P\}$ be the plane inversion 
defined by $I(Q)=P-\frac{1}{||Q-P||^2}(Q-P)$.    
Since $P\in \mathcal{NS}_{\theta, {\bf r}}$, 
the composed mapping ${\bf n}=I\circ f$ is well-defined \textcolor{black}{and of class $C^\infty$}.        
\textcolor{black}
{
The mapping ${\bf n}$ is exactly the dual curve of the family of lines 
$\left\{
{\bf r}(s)+R_\theta d{\bf r}_s(T_s(S^1))
\right\}_{s\in S^1}$ relative to the point $P\in \mathcal{NS}_{\theta, {\bf r}}$.
}
Let the boundary of convex hull of ${\bf n}(S^1)$ 
be denoted by $\partial{\rm conv}\left({\bf n}(S^1)\right)$.       
Then, by the construction, 
$\partial{\rm conv}\left({\bf n}(S^1)\right)$ intersect the half line 
$\{P+\lambda s\; |\; \lambda\in \mathbb{R}_+\}$ exactly at one point 
for any $s\in S^1$.     Thus, the composed image 
$I\left(\partial {\rm conv}\left(\textcolor{black}{\bf n}(S^1)\right)\right)$ intersect the half line 
$\{P+\lambda s\; |\; \lambda\in \mathbb{R}_+\}$ exactly at one point 
for any $s\in S^1$.       
Moreover, the intersecting points depend on $s$ continuously.    
Hence, by corresponding $s\in S^1$ to the distance between $P$ and 
the unique intersecting point $I\left(\partial {\rm conv}\left({\bf n}(S^1)\right)\right) \cap 
\{P+\lambda s\; |\; \lambda\in \mathbb{R}_+\}$, 
 we obtain  the well-defined continuous function $h_\theta: S^1\to \mathbb{R}_+$.     
Since $\textcolor{black}{\bf n}$ is of class $C^\infty$, it is easily seen that the obtained function $h_\theta$ satisfies the 
assumption of Theorem  6.3 in \cite{nishimurasakemi2}.     
Let $T_\theta: \mathbb{R}^2\to \mathbb{R}^2$ be the parallel translation given by 
$T_\theta(x,y)=(x, y)-P$.    
Then, by Theorem 6.3 of \cite{nishimurasakemi2}, it follows that    
\[
T_\theta( \overline{\mathcal{NS}_{\theta, {\bf r}}})=\mathcal{W}_{h_\theta}. 
\]
\hfill\qquad $\Box$ 
\begin{proposition}\label{proposition 5}
Let ${\bf r}: S^1\to \mathbb{R}^2$ be a $C^\infty$ immersion with its no-silhouette 
$\mathcal{NS}_{\bf r}$.    
Then, the map $\omega: [0, \theta_{\bf r})\to \mathcal{H}_{\mbox{\rm conv}}(\mathbb{R}^{2})$ 
defined by $\omega(\theta)=\overline{\mathcal{NS}_{\theta, {\bf r}}}$ is continuous, 
\end{proposition}
\underline{\it Proof of Proposition \ref{proposition 5}}\qquad 
Let $C^0(S^1, \mathbb{R}_+)$ be the set consisting of continuous function $S^1\to \mathbb{R}_+$.   
The set $C^0(S^1, \mathbb{R}_+)$ is a (non-complete) metric space with respect to the metric 
$d_{\mbox{\rm norm}}(h_1, h_2)=\max_{s\in S^1}|h_1(s)-h_2(s)|$.    
Let $\Gamma: [0, \theta_{\bf r})\to C^0(S^1, \mathbb{R}_+)$ 
(resp. $\Omega: C^0(S^1, \mathbb{R}_+)\to \mathcal{H}_{\rm conv}(\mathbb{R}^2)$) be the mapping 
defined by $\Gamma(\theta)=h_\theta$ 
(resp. $\Omega(h)=\mathcal{W}_h$), 
where $h_\theta$ is the continuous function defined in the proof of 
Proposition \ref{proposition 4}.       
Then, in order to show that $\omega$ is continuous, it is sufficient to show that  
both $\Gamma, \Omega$ are continuous.   
\par 
We first show that $\Gamma$ is continuous.   
Let $\widetilde{h}: S^1\to \mathbb{R}_+$ be the function defined by 
\[
\widetilde{h}(\cos \lambda, \sin\lambda)=
||P-I\left(\partial \mbox{\rm conv}\left(\textcolor{black}{\bf n}(S^1)\right)\right)
\cap \{P+t(\cos \lambda, \sin \lambda)\; |\; t\in \mathbb{R}_+\}||,    
\] 
where the set 
$I\left(\partial\mbox{\rm conv}\left(\textcolor{black}{\bf n}(S^1)\right)\right)
\cap \{P+t(\cos \lambda, \sin \lambda)\; |\; t\in \mathbb{R}_+\}$, 
which appeared in the proof of 
Proposition \ref{proposition 4}, 
is a one point set and it is regarded as a point.      
By obtaining the concrete expression of $\textcolor{black}{\bf n}$ given in the proof of Proposition \textcolor{black}{\ref{proposition 4}}, 
it is easily seen that $\textcolor{black}{\bf n}$ is smoothly depending on 
$\theta\in [0, \theta_{\bf r})$.     
Thus, $\widetilde{h}$ is continuously depending on $\theta\in [0, \theta_{\bf r})$.    
Since $I$ is a $C^\infty$ diffeomorphism of $\mathbb{R}^2-\{P\}$, it follows that 
$h_\theta$ is continuously depending on $\theta\in [0, \theta_{\bf r})$.    
Hence, $\Gamma$ is a continuous mapping.   
\par 
\medskip 
We next show that $\Omega$ is continuous.     
Let $\{h_i\}_{i=1,2,\ldots}\subset C^0(S^1, \mathbb{R}_+)$ be a convergent sequence to an element of $C^0(S^1, \mathbb{R}_+)$.    
Set $h=\lim_{i\to \infty}h_i$.     
We also set 
\[
W=
\left\{P\in \mathbb{R}^{2}\; \left| \; \exists P_i\in \mathcal{W}_{h_i}\; (i\in \mathbb{N}); \; 
\lim_{i\to \infty}P_i=P\right.\right\}.  
\]
Then, it is easily seen that $\mathbb{R}^2-W$ is an open set.  Thus, $W$ is a closed set.  
\par 
We show $\mathcal{W}_h=W$.     Let $P$ be an interior point of $\mathcal{W}_h$.   
Then, since $h=\lim_{i\to \infty}h_i$, $P$ must be an interior point of $\mathcal{W}_{h_i}$ 
for any sufficiently large $i$.    
Thus, $P$ is contained in $W$. 
Since both $\mathcal{W}_h$ and $W$ are closed, it follows that 
$\mathcal{W}_h\subset W$.    
Next, Let $Q$ be a point of $W$.    Suppose that 
$Q$ is not contained in $\mathcal{W}_h$.     
Then, there exists $s_0\in S^1$ such that 
$(Q\cdot s_0) > h(s_0)$, 
where $(Q\cdot s_0)$ 
stands for the scalar product of two vectors $Q, s_0\in \mathbb{R}^2$.    
Set $\varepsilon=(Q\cdot s_0)-h(s_0)>0$.    
Since $h=\lim_{i\to \infty}h_i$,  it follows that 
$(Q\cdot s_0)-h_i(s_0)>
\frac{\varepsilon}{2}$ for any sufficiently large $i$.   
This contradicts to the assumption that $Q\in W$.    Hence, we have that $W\subset \mathcal{W}_h$, and it follows that 
$\mathcal{W}_h=W$.    
\par 
The remaining part of the proof that $\Omega$ is continuous is to show the following:   
\[
\lim_{i\to \infty}d_H(W, \mathcal{W}_{h_i})=0.
\leqno{(2.2)}
\]       
In order to show (2.2), by the construction of $W$, it is sufficient to show that 
$\{\mathcal{W}_{h_i}\}_{i=1, 2, \ldots}$ is a Cauchy sequence of $\mathcal{H}(\mathbb{R}^2)$.     
Since $\{h_i\}_{i=1, 2, \ldots}$ is a Cauchy sequence of $C^0(S^1, \mathbb{R}_+)$, 
it is clear that $\{\mathcal{W}_{h_i}\}_{i=1, 2, \ldots}$ is a Cauchy sequence.    
Therefore, we have that 
$\lim_{i\to \infty}d_H(W, \mathcal{W}_{h_i})=0$ and it follows that 
$\Omega$ is continuous.        
\hfill\qquad $\Box$ 
\section{Proof of Theorem \ref{theorem 1}}\label{section 3} 
\textcolor{black}
{
By Proposition \ref{proposition 4}, there exists a parallel translation 
$T_\theta: \mathbb{R}^2\to \mathbb{R}^2$ such that 
$T_\theta\left(\overline{\mathcal{NS}_{\theta, {\bf r}}}\right)$ 
is a Wulff shape.      
In particular, $T_\theta\left(\overline{\mathcal{NS}_{\theta, {\bf r}}}\right)$ 
contains the origin as an interior point.   
}
Set 
$\widetilde{\bf r}=\alpha_N^{-1}\circ Id\textcolor{black}{\circ T_\theta}\circ {\bf r}$ and 
$\widetilde{\bf n}_\theta(s)=\cos\theta \widetilde{\bf n}(s)-\sin\theta \widetilde{\bf t}(s)$ for $s\in S^1$.  
We investigate the singularities of $\widetilde{\bf n}_\theta$.   
Let $U$ be an open arc of $S^1$.    
By using the arc-length parameter of $\widetilde{\bf r}|_U$, 
without loss of generality, from the first we may assume that $||\widetilde{\bf r}'(s)||=1$ for $s\in U$.   
Then, by Lemma \ref{lemma 2.1}, we have the following:   
\[
\widetilde{\bf n}'_\theta(s)=
-\kappa_g(s)\cos\theta\;\widetilde{\bf t}(s)+\sin\theta\; \widetilde{\bf r}(s) 
-\kappa_g(s) \sin\theta\; \widetilde{\bf n}(s).   
\] 
Since the angle $\theta$ satisfies $0<\theta<\theta_{\bf r}\le \frac{\pi}{2}$ in Theorem \ref{theorem 1}, 
it follows that 
$\sin\theta\ne 0$.     Therefore, $\widetilde{\bf n}_\theta$ is non-singular even at the point $s\in S^1$ 
such that $\kappa_g(s)=0$.    
\par 
Next, we show that $\widetilde{\bf n}_\theta(s)\cdot N>0$ for any $s\in S^1$.   
Let the dual of $\widetilde{\bf n}_\theta$ be denoted by $\widetilde{\bf r}_\theta$.    
Then,  it follows that $\widetilde{\bf r}_\theta$ is a Legendrian mapping and the following 
equality holds.   
\[
S^2_{N,+}\bigcap \left(S^2-\bigcup_{s\in S^1}GH_{\widetilde{\bf r}_\theta}\right)
=\alpha_N^{-1}\circ Id\circ
\mathcal{NS}_{\theta, {\bf r}}.  
\]
Since $\theta<\theta_{\bf r}$, by Lemma \ref{lemma 2.3}, we have that 
$\widetilde{\bf n}_\theta(s)\cdot N>0$ for any $s\in S^1$.   
Thus, the spherical convex hull of $\{\widetilde{\bf n}_\theta(s))\; |\; s\in S^1\}$ is well-defined.   
Since $\widetilde{\bf n}_\theta$ is non-singular, the boundary of  
$\mbox{\rm s-conv}(\{\widetilde{\bf n}_\theta(s))\; |\; s\in S^1\})$ is a submanifold of class $C^1$ 
(for instance see \cite{robertsonromerofuster, zakalyukin}).    
By the property (4) of $\Psi_N$, the boundary of 
$\Psi_N(\mbox{\rm s-conv}(\{\widetilde{\bf n}_\theta(s))\; |\; s\in S^1))$ is a submanifold of class $C^1$.   
It follows that the boundary of 
$Id^{-1}\circ\alpha_N\circ \Psi_N(\mbox{\rm s-conv}(\{\widetilde{\bf n}_\theta(s))\; |\; s\in S^1))$ is 
a submanifold of class $C^1$.   
\par 
On the other hand, 
\textcolor{black}
{
by constructions, it follows that 
}
$\textcolor{black}{T_\theta}(\overline{\mathcal{NS}_{\theta, {\bf r}}})$ is a Wulff shape $\mathcal{W}_h$ with the 
support function $h$ whose graph with respect to the polar coordinate expression 
is the boundary of 
$Id^{-1}\circ\alpha_N\circ \Psi_N(\mbox{\rm s-conv}(\{\widetilde{\bf n}_\theta(s))\; |\; s\in S^1))$. 
\par 
Therefore, the support function $h$ for the Wulff shape 
$\textcolor{black}{T_\theta}(\overline{\mathcal{NS}_{\theta, {\bf r}}})$ 
is of class $C^1$ and it follows that 
$\mathcal{W}_h$ is never a polygon by Proposition \ref{proposition 3}.    
\hfill\qquad $\Box$ 
\section{Proof of Theorem \ref{theorem 2}}\label{section 4} 
By Proposition \ref{proposition 4}, for any $i\in \mathbb{N}$ there exists a parallel translation 
$\textcolor{black}{T_{\theta_i}}: \mathbb{R}^2\to \mathbb{R}^2$ such that 
$\textcolor{black}{T_{\theta_i}}\left(\overline{\mathcal{NS}_{\theta_i, {\bf r}}}\right)$ 
is a Wulff shape $\mathcal{W}_{h_i}$ by an appropriate support function 
$\textcolor{black}{h_i}$.     
By Proposition \ref{dual Wulff shape},  for any $i\in \mathbb{N}$ the set 
$Id^{-1}\circ \alpha_N\left(\left(\alpha_N^{-1}\circ Id(\mathcal{W}_{h_i})\right)^\circ\right)$ is a Wulff shape too.   
Thus, by Proposition \ref{taylor}, it follows that both 
$\alpha_N^{-1}\circ Id(\mathcal{W}_{h_i})$ and  
$\left(\alpha_N^{-1}\circ Id(\mathcal{W}_{h_i})\right)^\circ$ belong to 
$\mathcal{H}(S^2)$ for any $i\in \mathbb{N}$.     
Moreover, \textcolor{black}{by} Proposition \ref{proposition 5}, 
we may assume that $\{T_{\theta_i}\}_{i=1, 2, \ldots}$ is a Cauchy sequence.   
Thus, \textcolor{black}{we may assume that} 
both $\left\{\alpha_N^{-1}\circ Id(\mathcal{W}_{h_i})\right\}_{i=1, 2, \ldots}$ and  
$\left\{\left(\alpha_N^{-1}\circ Id(\mathcal{W}_{h_i})\right)^\circ\right\}_{i=1, 2, \ldots}$ are Cauchy sequences.     
\par 
By Proposition \ref{proposition 2}, $\lim_{i\to \infty}\overline{\mathcal{NS}_{\theta_i, {\bf r}}}$ is a 
point or segment.    Suppose that it is a segment.      
Let $P_1, P_2\in S^2$ be two boundary points of this segment.   
Then, by Proposition \ref{continuity} and Lemma \ref{Maehara's lemma}, we have the following:    
\[
\lim_{i\to \infty}\left(\alpha_N^{-1}\circ Id(\mathcal{W}_{h_i})\right)^\circ 
= 
H(P_1)\cap H(P_2).   
\] 
\par 
Let $\widetilde{\bf n}_{\theta_{\bf r}}: S^1\to S^2$ be the $C^\infty$ mapping defined by 
$\widetilde{\bf n}_{\theta_{\bf r}}(s)=\cos \theta_{\bf r}\widetilde{\bf n}(s)-\sin \theta_{\bf r}\widetilde{\bf t}(s)$ for any $s\in S^1$, 
where $\widetilde{\bf n}$ and $\widetilde{\bf t}$ are the same $C^\infty$ mapping as in the proof of Theorem \ref{theorem 1}.   
Then, notice that 
$\widetilde{\bf n}_{\theta_{\bf r}}(S^1)\subset \textcolor{black}{H}(P_1)\cap H(P_2)$.   
For any $j$ $(j=1,2)$, we let the set $\{Q\in S^2\; |\; P_j\cdot Q=0\}$ be denoted by $\partial H(P_j)$.     
Then, the intersection $\partial H(P_1)\cap \textcolor{black}{\partial} H(P_2)$ consists of two antipodal points $Q_1, Q_2$.    
By Lemma \ref{lemma 2.3} and Proposition \ref{proposition 1}, there exists $s_1, s_2\in S^1$ $(s_1\ne s_2)$ such that 
$\widetilde{\bf n}_{\theta_{\bf r}}(s_1)=Q_1$, $\widetilde{\bf n}_{\theta_{\bf r}}(s_2)=Q_2$.      
\par 
On the other hand, since $0\le \theta_{\bf r}\le \frac{\pi}{2}$, similarly as in the proof of Theorem \ref{theorem 1}, 
it follows that $\widetilde{\bf n}_{\theta_{\bf r}}$ is non-singular.    Thus, we have a contradiction.   
\hfill\qquad $\Box$  
\section{Proof of Theorem \ref{theorem 3}}\label{section 5}
For any $\theta$ $(0\le \theta<\theta_{\bf r})$ and any $s\in S^1$, set 
\[
\ell_{\theta, s}={\bf r}(s)+R_\theta\left(d{\bf r}_s(T_sS^1)\right).   
\]
Let $f_{\theta, s}(x,y)$ be the affine function which define $\ell_{\theta, s}$.    
Set 
\[
H_{\theta, s}^+=\{(x,y)\in \mathbb{R}^2\; |\; f_{\theta, s}(x,y)>0\}, \quad 
H_{\theta, s}^-=\{(x,y)\in \mathbb{R}^2\; |\; f_{\theta, s}(x,y)<0\}.   
\]
Then, since $\textcolor{black}{\overline{\mathcal{NS}_{\theta, {\bf r}}}}$ 
is a convex body for any $\theta$ 
$(0\le \theta <\theta_{\bf r})$, it follows that one of 
$\mathcal{NS}_{\theta, {\bf r}}=\cap_{s\in S^1}H_{\theta, s}^+$ or 
$\mathcal{NS}_{\theta, {\bf r}}=\cap_{s\in S^1}H_{\theta, s}^-$ \textcolor{black}{holds}.     
By Proposition \ref{proposition 4}, we may assume that the following holds 
for any $\theta$ $(0\le \theta <\theta_{\bf r})$.    
\[
\mathcal{NS}_{\theta, {\bf r}}=\bigcap_{s\in S^1}H_{\theta, s}^+.   
\]
Since ${\bf r}$ does not have inflection points, it follows that 
$\mathcal{NS}_{0, {\bf r}}$ contains $\mathcal{NS}_{\theta, {\bf r}}$ for any 
$\theta$ $(0\le \theta <\theta_{\bf r})$.    Thus, for any $\theta$ $(0\le \theta <\theta_{\bf r})$, 
we have the following:   
\begin{eqnarray*}
\mathcal{NS}_{\theta, {\bf r}} & = & 
\mathcal{NS}_{\theta, {\bf r}}\cap \mathcal{NS}_{0, {\bf r}} \\ 
{ } & = & 
\left(\bigcap_{s\in S^1}H_{\theta, s}^+\right)\bigcap \mathcal{NS}_{0, {\bf r}} \\ 
{ } & = & 
\bigcap_{s\in S^1}\left(H_{\theta, s}^+\bigcap \mathcal{NS}_{0, {\bf r}}\right).    
\end{eqnarray*}
Since ${\bf r}$ does not have inflection points, we have that 
$H_{\theta_1, s}^+\cap \mathcal{NS}_{0, {\bf r}}$ contains 
$H_{\theta_2, s}^+\cap \mathcal{NS}_{0, {\bf r}}$ for any 
two $\theta_1, \theta_2\in [0,\theta_{\bf r})$ satisfying 
$0\le \theta_1<\theta_2<\theta_{\bf r}$.       It follows that 
$\mathcal{NS}_{\theta_1, {\bf r}}\supset \mathcal{NS}_{\theta_2,{\bf r}}$ 
if $0\le\theta_1<\theta_2<\theta_{\bf r}$.   
\hfill\qquad $\Box$ 
\section*{Acknowledgements}
\textcolor{black}
{
The authors would like to thank the referee for careful reading of the first draft of their paper and giving some comments.       
}
\par 
T.~Nishimura is partially supported by 
JSPS KAKENHI Grant Number 26610035. 


\end{document}